\documentclass[a4paper,12pt, leqno ]{article}

\title{Some remarks on the geometry of Kropina spaces. }
\author{By Ryozo Yoshikawa and Sorin V. Sabau}
\date{}
\newlength{\topdummy}
\makeatletter
 
  \@addtoreset{equation}{section}
 \makeatother

\usepackage{latexsym}
\usepackage{amssymb}
\usepackage{enumerate}
\usepackage{epic,eepic}
\usepackage{color}
\usepackage{float}

\newcommand{\R}{\ensuremath{\mathbb{R}}}

\begin{document}

\newtheorem{definition}{Definition}[section]
\newtheorem{theorem}{Theorem}[section]
\newtheorem{proposition}[theorem]{Proposition}
\newtheorem{corollary}[theorem]{Corollary}
\newtheorem{lemma}[theorem]{Lemma}
\newtheorem{remark}[theorem]{Remark}
\newtheorem{example}{Example}[section]

\maketitle
\begin{abstract}

Using the navigation data $(h, W)$ of a Kropina space $(M, \alpha^2/\beta)$, 
we characterize
weakly-Berwald Kropina spaces and Berwald  Kropina spaces by means of the Killing vector field $W$ and the parallel vector field $W$, respectively.  

Moreover, the local $1$-parameter group of Finslerian local isometries on $(M, \alpha^2/\beta)$ coincides with the local $1$-parameter group Riemannian of local isometries on $(M, h)$.
\end{abstract}

\section{Introduction}

Finsler metrics generalize Riemannian metrics. The most natural Finsler structures are those obtained by deformations of Riemannian metrics. Randers metrics are the most famous metrics of this type because of their relation with the Zermelo's navigation problem. 
On the other hand, recently it was shown that Kropina metrics also give solutions to the Zermelo's navigation problem (\cite{YS} and references herein). This suggests that Kropina spaces are Finsler spaces with rich geometrical properties and that their geometry worth more investigations.

In the present paper we present three remarks on the geometry of Kropina spaces.

The first remark concerns the relation between weakly-Berwald, Berwald and Kropina spaces. We prove that the set of weakly-Berwald Kropina metrics coincide with the set of strong Kropina metrics defined in the paper (Theorem 3.2). Moreover, a Kropina space is a Berwald one if and only if the wind vector field $W$ is parallel with respect to the Riemannian metric $h$ (Theorem 3.3). Nevertheless, in the Kropina case,  Landsberg and Berwald spaces coincide. 

The second remark is about the characterization of Kropina metrics of $p$-scalar curvature, i.e. Kropina metrics whose scalar flag curvature  is a function of position only $K(x,y)=K(x)$. We give a characterization of these spaces in terms of navigation data $(h,W)$ (Theorem 4.1). Moreover, we show that Kropina spaces of $p$-scalar curvature must be Berwald  if and only if $K(x)=0$ (Theorem 4.4). 

The third remark is on the isometry group of a Kropina metric. In general, the isometry group of a Finsler structure is not easy to be determined. However, in the case of strong Kropina metrics we prove that its local isometry group actually coincides with the local isometry group of the Riemannian space $(M,h)$. 

\section{The navigation data of Kropina spaces.}
A Kropina metric $F=\alpha^2/\beta$, where $\alpha=\sqrt{a_{ij}(x)y^iy^j}$ and $\beta=b_i(x)y^i$, on an $n$-dimensional differential manifold $M$ is characterized by a new Riemannian metric $h=\sqrt{h_{ij}(x)y^iy^j}$ on $M$ and a vector field $W=W^i(\partial/\partial x^i)$ of constant length 1 with respect to $h$, where
\begin{eqnarray}\label{equation 1.1}
    h_{ij}:=e^{\kappa(x)}a_{ij}, \hspace{0.1in} 2W_i:=e^{\kappa(x)}b_i, \hspace{0.1in} e^{\kappa(x)}b^2=4 
\end{eqnarray}
for a function $\kappa(x)$ of $(x^i)$ alone.
In the above equations, we used $W_i(x):=h_{ij}(x)W^j(x)$ and $b^2:=a^{ij}b_ib_j$, where $(a^{ij})$ is the inverse matrix of $(a_{ij})$ (\cite{YO1}, \cite{YO2}, \cite{YO3}, \cite{YS}).

The pair $(h, W)$ is called the {\it navigation data} of the Kropina metric $F=\alpha^2/\beta$ or the Kropina space $(M, F)$. 
The following theorem is very important:
\begin{theorem}(\cite{YO2}, \cite{YO3})\label{Theorem 1.1}
 A Kropina space $(M, F=\alpha^2/\beta)$ whose navigation data is $(h, W)$ is  of constant flag curvature $K$ if and only if the  following two conditions holds:
\begin{enumerate}[(a)]
\item $W$ is a unit Killing vector field,
\item the Riemannian space $(M, h)$ is of constant sectional curvature $K$.
\end{enumerate}
\end{theorem}

We have 
\begin{definition}
Let  $(M, F=\alpha^2/\beta)$ be a Kropina space  whose navigation date is $(h, W)$. 
If $W$ is a unit Killing vector field, we call $(M, F)$ a {\it strong Kropina   space}.
\end{definition}

\section{ Berwald spaces and weakly-Berwald spaces}

In this section, we will consider Kropina spaces which are weakly-Berwald spaces or Berwald spaces.

In general,  we denote the coefficients  of the geodesic spray of a Finsler space by $G^i$  and put
\begin{eqnarray*}
  {G^i}_j:=\frac{\partial G^i}{\partial y^j}, \hspace{0.2in} {{G_j}^i}_k:=\frac{\partial {G^i}_j}{\partial y^k}, \hspace{0.2in} {{G_j}^i}_{kl}:=\frac{\partial {{G_j}^i}_k}{\partial y^l}.
\end{eqnarray*}
Berwald spaces are defined as follows:
\begin{definition}
A Finsler space $(M, F)$ is called a Bewald one if $G^i$ satisfy the conditions  ${{G_j}^i}_{kl}=0$, that is, the coefficients  ${{G_j}^i}_k$ of the Berwald connection of $(M, F)$ are functions of the position $(x^i)$ alone.
\end{definition}

The concept of weakly-Berwald spaces is a generalization of that of Berwald spaces.

\begin{definition}
If a Finsler space satisfies the condition $G_{ij}:={{G_i}^r}_{jr}=0$, we call it a {\it weakly-Berwald one}.
\end{definition}

The first author and Sandor B\'acs\'o considered the relation between the concept of weakly-Berwald spaces and some other concepts in \cite{BY}. 
Furthermore, the same author, Katsumi Okubo and Makoto Matsumoto have obtained the necessary and sufficient conditions for a Kropina space to be a weakly-Berwald one or a Berwal one in \cite{YOM}.

\begin{theorem}(\cite{YOM})\label{Theorem 2.1}
Let $(M, \alpha^2/\beta)$ be a Kropina space, where $\alpha=\sqrt{a_{ij}(x)y^iy^j}$ and $\beta=b_i(x)y^i$.
The necessary and sufficient conditions for a  Kropina space to be a weakly-Berwald one or a Berwald one are as follows respectively:\\
$(wB)$  \hspace{0.3in}   $r_{ij}=c(x)a_{ij}$ for a function $c(x)$ of the position $(x^i)$ alone,\\
$(B)$  \hspace{0.4in}  $r_{ij}=c(x)a_{ij}$ for a function $c(x)$ of the position $(x^i)$ alone and $s_jb_i-s_ib_j=b^2s_{ij}$.

In the above equations, we used   $r_{ij}:=(b_{i;j}+b_{j;i})/2$,  $s_{ij}:=(b_{i;j}-b_{j;i})/2$, $s_i=a^{rs}b_rs_{si}$, 
where  the symbol $(;)$ denotes the covariant derivative with respect to the Riemannian  metric $\alpha$ and $(a^{rs})$ is the inverse matrix of $(a_{rs})$.
\end{theorem}

First, we suppose that a Kropina space $(M, F=\alpha^2/\beta)$ is a weakly-Berwald one. Then, from Theorem \ref{Theorem 2.1}, we have  $r_{ij}=c(x)a_{ij}$ for a function $c(x)$ of the position $(x^i)$ alone.
It   can be changed as follows:
\begin{eqnarray}\label{equation 2.1}
r_{ij}=c(x)e^{-\kappa}h_{ij}.
\end{eqnarray}

Now, from Section 2 in \cite{YO2}  we have an equation
\begin{eqnarray}\label{equation 2.2}
  r_{ij}=2e^{-\kappa}\bigg(\texttt{R}_{ij}-\frac{1}{2}W_r\overline{\kappa}^rh_{ij}\bigg),
\end{eqnarray}
where
$\texttt{R}_{ij}:=(W_{i||j}+W_{j||i})/2$,  $\overline{\kappa}^r:=h^{rs}(\partial \kappa/\partial x^s)$,  the symbol $(||)$ denotes the covariant derivative with respect to the metric $h$ and $(h^{rs})$ is the inverse matrix of $(h_{rs})$.

Substituting (\ref{equation 2.1}) to (\ref{equation 2.2}), we get
\begin{eqnarray}\label{equation 2.3}
 \texttt{R}_{ij}=\frac{1}{2} \bigg(c(x) +W_r\overline{\kappa}^r\bigg)h_{ij}.
\end{eqnarray}
Since $W$ is a unit vector field on $(M, h)$, we have $|W|^2=W^iW^jh_{ij}=1$ and $W_{i||k}W^i=0$.
So, transvecting (\ref{equation 2.3}) by $W^iW^j$, we get
     $0=\texttt{R}_{ij}W^iW^j=\frac{1}{2} \bigg(c(x) +W_r\overline{\kappa}^r\bigg)|W|^2$,
that is, 
              $c(x) +W_r\overline{\kappa}^r=0$.
Hence, the equation (\ref{equation 2.3}) reduces to 
\begin{eqnarray}\label{equation 2.4}
\texttt{R}_{ij}=0,
\end{eqnarray}
that is, the vector field $W$ is Killing.

Conversely, suppose that the vector field $W$ is a Killing one.
Substituting (\ref{equation 2.4}) to (\ref{equation 2.2}), we get
   $r_{ij}=-e^{-\kappa}W_r\overline{\kappa}^rh_{ij}=-W_r\overline{\kappa}^ra_{ij}$.
From Theorem \ref {Theorem 2.1}, the Kropina space $(M, F=\alpha^2/\beta)$ is a weakly-Berwald one.

So, we obtain
\begin{theorem}\label{Theorem 2.2}
Let $(M, F=\alpha^2/\beta)$ be a Kropina  space. Let $(h, W)$ be the navigation data of $F=\alpha^2/\beta$ defined by (\ref{equation 1.1}).
Then the Kropina space  $(M, F=\alpha^2/\beta)$ is a weakly-Bewald one if and only if it is a strong Kropina one.
\end{theorem}

In fact, if we suppose that the vector field $W$ is Killing, from Theorem 5 in \cite{YO2} we have 
\[ 2G^i=^h{{\gamma_0}^i}_0-2F{\texttt{S}^i}_0,\]
where $^h{{\gamma_0}^i}_0$ are the Christoffel's symbols of the Riemannian metric $h$, $\texttt{S}_{ij}:=(W_{i||j}-W_{j||i})/2$, ${\texttt{S}^i}_j:=h^{ir}\texttt{S}_{rj}$ and $(h^{ij})$ is the inverse matrix of $(h_{ij})$.
From the above equation, we get
\begin{eqnarray}\label{equation 2.5}
2{{G_j}^i}_k &=&2^h{{\gamma_j}^i}_k-\frac{2h_{jk}}{W_0}{\texttt{S}^i}_0   +\frac{2h_{0j}W_k}{(W_0)^2}{\texttt{S}^i}_0    
                                              +\frac{2h_{0k}W_j}{(W_0)^2}{\texttt{S}^i}_0 -\frac{2h_{00}W_jW_k}{(W_0)^3}{\texttt{S}^i}_0             \nonumber\\
              &&-\frac{2h_{0j}}{W_0}{\texttt{S}^i}_k-\frac{2h_{0k}}{W_0}{\texttt{S}^i}_j   +\frac{h_{00}W_j}{(W_0)^2}{\texttt{S}^i}_k+\frac{h_{00}W_k}{(W_0)^2}{\texttt{S}^i}_j. \nonumber
\end{eqnarray}

Putting $i=k=r$, we obtain
     ${{G_j}^r}_r=^h{{\gamma_j}^r}_r$.
Since $^h{{\gamma_j}^r}_r$ is a function of  $(x^i)$ alone, it follows that the Kropina space is a weakly-Berwald space.\\

Next, suppose that a Kropina space $(M, F=\alpha^2/\beta)$ is a Berwald one.  
We consider the two conditions in $(B)$ of Theorem \ref{Theorem 2.1}.
By straightforward calculations,   it follows that the second condition  $s_jb_i-s_ib_j=b^2s_{ij}$  can be changed as follows:
\begin{eqnarray}\label{equation 2.6}
  \texttt{S}_{ij}=W_i\texttt{S}_j-W_j\texttt{S}_i.
\end{eqnarray}

From Theorem \ref{Theorem 2.2},   the first condition  means  that the unit vector field $W$ is Killing, that is, $\texttt{R}_{ij}=0$.
Hence, we get $\texttt{S}_i=W^r\texttt{S}_{ri}=   W^rW_{r||i}=0$.
Substituting it to (\ref{equation 2.6}), we get
   $ \texttt{S}_{ij}=0$.
From it and $\texttt{R}_{ij}=0$, we get $W_{i||j}=0$, that is, $W$ is a parallel vector field on the Riemannian space $(M, h)$.

Conversely, suppose that the  unit vector field $W$ is a parallel vector field on the Riemannian space $(M, h)$.
Then, the equations  $\texttt{R}_{ij}=0$ and (\ref{equation 2.6}) hold.
 Hence, from Theorem \ref{Theorem 2.1} the Kropina space is a Berwald space. 
In fact, if we suppose that $W$ is parallel,  we get  ${{G_j}^i}_k=^h{{\gamma_j}^i}_k$ because of $\texttt{S}_{ij}=0$.

Therefore, we get

\begin{theorem}\label{Theorem 2.3}
Let $(M, F=\alpha^2/\beta)$ be a Kropina space, where $\alpha=\sqrt{a_{ij}(x)y^iy^j}$ and $\beta=b_i(x)y^i$.
And let $(h, W)$ be the navigation data of $F$ defined by (\ref{equation 1.1}).

 Then, the Kropina space $(M, \alpha^2/\beta)$ is a Berwald space if and only if $W$ is parallel with respect to  the Levi-Civita connection  of the metric $h$. 

In this case,  the coefficients ${{G_j}^i}_k$ of the Berwald connection coincides with the coefficient $^h{{\gamma_j}^i}_k$ of the Levi-Civita connection of $h$.
In other words, the geodesics of a Kropina space $(M, F=\alpha^2/\beta)$ which is a Berwald one coincide with those of a Riemannian space $(M, h)$.
\end{theorem}

\begin{remark}
We consider a Finsler space $(M, F=\alpha^2/\beta)$ whose navigation data is $(h, W)$.
If it is a Berwald space, the vector field $W$ satisfies the condition $W_{i||j}=0$. Hence, it is a strong Kropina space. That is, the set of Kropina spaces which are Berwald spaces are contained in the set of the strong Kropina spaces. 
\end{remark}

\section{A Kropina space of p-scalar flag curvature}
We consider a Kropina space of scalar flag curvature $K=K(x)$ which is a scalar function of the position $(x^i)$ alone. 
So, we have
\begin{definition}
Let a (conic) Finsler space $(M, F)$ be of scalar flag curvature $K=K(x, y)$.
If $K=K(x, y)$ is independent of $(y^i)$ and is a function of the position $(x^i)$ alone, the (conic)  Finsler space $(M, F)$ is called to be  {\it of $p$-scalar flag curvature}.
\end{definition}

Observing the proof of Theorem 4 in \cite{YO2}, it follows that even if we exchange "of constant flag curvature $K$" with " of $p$-scalar flag curvature $K=K(x)$",  the Theorem  4 holds good.
 Namely, we obtain
\begin{theorem}\label{Theorem 3.1}
 A Kropina space $(M, F=\alpha^2/\beta)$ whose navigation data is $(h, W)$ is of $p$-scalar flag curvature $K=K(x)$  if and only if the  following two conditions holds:
\begin{enumerate}[(a)]
\item The unit vector field $W$ is   Killing,
\item the Riemannian space $(M, h)$ is of scalar sectional curvature $K=K(x)$.
\end{enumerate}

Furthermore, in this case  $K=K(x)\ge 0$.
\end{theorem}

In the above theorem, we must remark that if the dimension of $M$ is more than 2, from Shur's Lemma it follows that $K(x)$ is constant.

\begin{remark}

Pay attention to Berger's theorem (\cite{BN}, \cite{Ber}, \cite{V}):\\

   {\it Every Killing vector field on a compact even-dimensional Riemannian space $(M, g)$ with positive sectional curvature vanishes.}\\

Then, it follows that there  exists no  Kropina space of $p$-scalar flag curvature on a two-dimensional compact Riemannian space.
  
\end{remark}
We shall prove that if  a Kropina space $(M, F=\alpha^2/\beta)$ is of $p$-scalar flag curvature $K=K(x)$, we have $K=K(x)\ge 0$.

Since the vector field $W$ is of  unit length,   we have 
\begin{eqnarray}\label{equation 3.2}
W_{r||i}W^r=0.
\end{eqnarray}
 From the equation (\ref{equation 3.2}), we get
\begin{eqnarray}\label{equation 3.3}
W_{r||i}{W^r}_{||j}+W_{r||i||j}W^r=0.
\end{eqnarray}

We have 
\begin{lemma}(\cite{YO2})
Let $(M, h)$ be a Riemannian space.
For a unit Killing vector field $W=W^i(\partial/\partial x^i)$ on $M$, the equality
\begin{eqnarray}\label{equation 3.4}
W_{i||j||k}={W_r}^h{{R_k}^r}_{ij}
\end{eqnarray}
holds good. 
In the above equation, the symbol $(_{||})$ denotes the covariant derivative with respect to the metric $h$ and $^h{{R_k}^r}_{ij}$ denotes the curvature tensor of the Riemannian space $(M, h)$.
\end{lemma}

Since the Riemannian space $(M, h)$ is of scalar sectional curvature $K=K(x)$, we have
\begin{eqnarray}\label{equation 3.5}
^h{{R_k}^r}_{ji}=K(x)(h_{kj}{\delta^r}_i-h_{ki}{\delta^r}_j).
\end{eqnarray}

From (\ref{equation 3.3}),  (\ref{equation 3.4})  and (\ref{equation 3.5}), we get
\begin{eqnarray*}
W_{r||i}{W^r}_{||j}=-W_{r||i||j}W^r=-{W_s}^h{{R_j}^s}_{ri}W^r=-K(x)W_s(h_{jr}{\delta^s}_i-h_{ji}{\delta^s}_r)W^r=K(x)(h_{ji}-W_jW_i),
\end{eqnarray*}
that is,
\begin{eqnarray}\label{equation 3.6}
    h_{rs}{W^s}_{||i}{W^r}_{||j}=K(x)(h_{ij}-W_iW_j).
\end{eqnarray}
Transvecting (\ref{equation 3.6}) by $y^iy^j$, we get
\begin{eqnarray*}
    h_{rs}{W^s}_{||0}{W^r}_{||0}=K(x)\{h_{00}-(W_0)^2\}.
\end{eqnarray*}
If we suppose that the equation ${W^s}_{||0}=0$ holds, we have  $K=K(x)=0$ because of $h_{00}-(W_0)^2\ne 0$.
On the other hand, if we suppose that ${W^s}_{||0}\ne 0$,  we have  $h_{rs}{W^s}_{||0}{W^r}_{||0}> 0$   since $(h_{rs})$ is positive definite.
From Schwarz's inequality, we have $h_{00}-(W_0)^2> 0$.
Hence, we obtain $K=K(x)> 0$.\\

From the above proof, we obtain
\begin{theorem}
Let  $(M, F=\alpha^2/\beta)$ be a Kropina space of $p$-scalar flag curvature $K=K(x)$ and  $(h, W)$  be the navigation data of $F$.
Then  $(M, F=\alpha^2/\beta)$ is a Berwald space, that is, ${W^i}_{||j}=0$ if and only if $K=K(x)=0$.  
\end{theorem}

\section{Finslerian isometries}
Let $(M, F=\alpha^2/\beta)$ be a Kropina space and $(h, W)$ be its navigation data. 
Let $\{\varphi_t\}$ be a local $1$-parameter group of local transformations on $M$. 
In this section, we consider the relation between the following two statements:

(1) \hspace{0.2in}  $\varphi_t$$(t\in I_\epsilon)$ are Riemannian local isometries on $(M, h)$,

(2) \hspace{0.2in} $\varphi_t$$(t\in I_\epsilon)$ are Finslerian  local isometries on $(M, F=\alpha^2/\beta)$.

\subsection{ Killing vector fields}
First, we define a {\it 1-parameter group of transformation} of $M$.
\begin{definition}\label{Definition 4.1}(\cite{Kl})
Let $M$ be a differential manifold.
A 1-parameter group of transformations is a differential mapping
\begin{eqnarray*}
        \varphi : \R \times M \longrightarrow M ; \hspace{0.1in} (t, p) \longmapsto t\cdot p \equiv \varphi_t(p)
\end{eqnarray*}
such that
\begin{enumerate}[(i)]
\item  $t \cdot (t'\cdot p)=(t+t')\cdot p$;  $0\cdot p=p$,
\item  $\varphi_t : M \longrightarrow M$ ;  $p \longmapsto t\cdot p$ is a transformation for every $t\in \R$.
\end{enumerate}
$\varphi$ determines a vector field $X_{\varphi}$ on $M$ by
\begin{eqnarray*}
X_\varphi(p):=\frac{d}{dt}\varphi_t(p)\bigg|_{t=0}.
\end{eqnarray*}
\end{definition} 

\begin{definition}(\cite{Kl})
Let  $M$ be a Riemannian manifold $M=(M, g)$ and $\varphi_t$ defined in Definition \ref{Definition 4.1} be an isometry for all $t$. 
Then $\{\varphi_t\}_{t\in \R}$ is called 1-parameter group of isometries. The associated vector field $X_\varphi$ is called Killing vector field. 
\end{definition}

On the converse of  Definition \ref{Definition 4.1} we have the following theorem :

\begin{theorem}(\cite{KN})\label{Theorem 4.1}
Let $X$ be a vector field on a manifold $M$.
For each point $p$ of $M$, there exist a neighborhood $U$ of $p$, a positive number $\epsilon$ and a local $1$-parameter group of local transformations $\varphi_t : U \longrightarrow M$, $t\in I_\epsilon=(-\epsilon, \epsilon)$ which induces the given $X$.
\end{theorem}

\begin{definition}(\cite{KN})
Let $M$ be a differential manifold and $X$ be a vector field on $M$.
If there exists a global $1$-parameter group of transformations of $M$ which induces  $X$, we say that $X$ is complete.
\end{definition}

We must remark that if $\varphi_t(p)$ is defined on $I_\epsilon \times M$ for some $\epsilon$, the $X$ is complete.

From the theory of the ordinary differential equations, we can give the following remark:
\begin{remark}\label{Remark 4.2}(\cite{KN})
  If two local $1$-parameter group of local transformations $\varphi_t$ and $\psi_t$ defined on $I_\epsilon \times U$ induce the same vector field on $U$, they coincide on $U$.
\end{remark} 

Let $X$ be a globally defined vector field on $M$. From Theorem \ref{Theorem 4.1},  for any point $p\in M$ there exist a neighborhood $U_p$ of $p$, a positive number $\epsilon_p$ and a local $1$-parameter group of local transformations $\varphi_t : U_p \longrightarrow M$, $t\in I_{\epsilon_p}$ which induces the given $X$. 
For another point $q\in M$ which is not $p$, there exist a neighborhood $U_q$ of $q$, a positive number $\epsilon_q$ and a local $1$-parameter group of local transformations $\psi_t : U_q \longrightarrow M$, $t\in I_{\epsilon_q}$ which induces the same $X$. Suppose that $U_p \cap U_q \ne \phi$. Then  since $\varphi_t$ and $\psi_t$ generate the same vector field $X$ on $U_p \cap U_q$, it follows that they are coincide on  $U_p \cap U_q$ from Remark \ref{Remark 4.2}. 

But even if we take  $t$ such that $|t|$ is  sufficiently small,  we cannot necessarily extend a local transformation $\varphi_t$  to a global transformation on $M$.  

We have the following proposition:
\begin{proposition}(\cite{KN})\label{Prop. compact}
On a compact manifold $M$, every vector field $X$ is complete.
\end{proposition}

Let $(M, F=\alpha^2/\beta)$ be a strong Kropina space whose navigation data is $(h, W)$. Since $W$ is a Killing vector field on $M$, for any point $p\in M$ there exists a local $1$-parameter group $\{\varphi_t\}_{t\in I_\epsilon}$ of   local isometries on $M$ which generate $W$. 
Since $\varphi_t$ is a local isometry,   we have
\begin{eqnarray}\label{equation 4.1}
   ( (\varphi_{-t})^*h_{\varphi_t(p)} )(Y_p, Z_p)=h_{\varphi_t(p)}((\varphi_t)_*(Y_p), (\varphi_t)_*(Z_p))=h_p(Y_p, Z_p)
\end{eqnarray}
for any vector fields $Y$ and $Z$ on $M$.
Using (\ref{equation  4.1}), we have
\begin{eqnarray*}
   (L_Wh)_p(Y_p, Z_p) =\lim_{t\longrightarrow  0}\frac{(\varphi_{-t})^*h_{\varphi_t(p)}-h_p}{-t}(Y_p, Z_p) =0 
\end{eqnarray*}  
for any vector fields $Y$ and $Z$ on $M$.

Conversely, if we suppose that $\{\varphi_t\}_{t\in I_\epsilon}$ is a local $1$-parameter group of local transformations on $M$ and that  the equation $(L_Wh)_p=0$ holds for every point $p\in M$, it follows that $\varphi_t$ is a local isometry of $M$. 

Then we have

\begin{proposition}\label{Proposition 4.4}
Let $\{\varphi_t\}_{t\in I_\epsilon}$ be a local  $1$-parameter group of local transformations on a Riemannian manifold $(M, h)$ and $W$ be a vector field on $M$ generated by $\{\varphi_t\}_{t\in I_\epsilon}$.  Then $W$ is a Killing vector field if and only if the equation $L_Wh=0$ holds. 
\end{proposition}

\subsection{Finslerian isometries}
Let $(M, F=\alpha^2/\beta)$ be a  Kropina space whose navigation data is $(h, W)$.
From Theorem \ref{Theorem 4.1}, for any point $p\in M$ there exist a neighborhood $U$ of $p$, a positive number $\epsilon$ and a local $1$-parameter group of local transformations $\varphi_t : U \longrightarrow M$ $(t\in I_\epsilon)$ which induce $W$.

Let $(x^i)$ be a local coordinate system of a  neighborhood $U$ of $p$.
 Suppose that $|t|$ is sufficiently small such that $\varphi_t(p)\in U$.  
We denote $p=(x^i)$ and $\varphi_t(p)=(\overline{x}^i)$.
Since $\{\varphi_t\}_{t\in I_\epsilon}$ generates a vector field $W=W^i(\partial/\partial x^i)$ on   a  neighborhood $U$ of $p$, we have 
  \[\overline{x}^i =x^i+tW^i+O(t^2).\]
From the above equation,  we get
\begin{eqnarray}\label{equation 4.8}
  \overline{y}^i &=&y^i+t\frac{\partial W^i}{\partial x^j}y^j+O(t^2). 
\end{eqnarray}

We define isometries  and local isometries on a Finsler space $(M, F)$ as follows:

\begin{definition}
Let $\phi$ be a transformation on an $n$-dimensional differential manifold $M$.
Then we call $\phi$  a Finslerian isometry on a Finsler manifold $(M, F)$ if $F(\phi(x), \phi(y))=F(x, y)$ for all $(x, y)\in TM$.
\end{definition}

\begin{definition}
Let $U$ be an open set of $n$-dimensional differential manifold $M$.
Let $\phi$ be a local transformation  from  $U$ to $M$.
Then we call $\phi$  a Finslerian local isometry on $U$ if $F(\phi(x), \phi(y))=F(x, y)$ for all $(x, y)\in TU$.
\end{definition}

We have 
\begin{theorem}\label{Theorem 4.3}
Let $(M, F)$ be a Finsler manifold and $\phi$ be a Finslerian local  isometry on an open set $V$ of $M$.
Suppose that a curve $\wp(t)$ is a geodesic on $V$, then the image $\phi(\wp(t))$  of $\wp(t)$ by $\phi$ is a geodesic on $(M, F)$.
\end{theorem}

Let $(M, F=\alpha^2/\beta)$ be a Kropina space and $(h, W)$ be a navigation data of $F$.
Suppose that $\{\varphi_t\}_{t\in I_\epsilon}$ is a local $1$-parameter group of  local isometries   on $M$.

We consider the necessary and sufficient conditions for  $\{\varphi_t\}_{\epsilon>t>0}$  to be a local $1$-parameter group of  Finslerian local transformations on  $(M, F=\alpha^2/\beta)$.
Then we have to only obtain the conditions for a local transformation $\varphi_t$ $(\epsilon>t>0)$  to be a Finslerian local isometry.

Let $p\in M$ and $U$ be a neighborhood of $p$.  
Suppose that $\varphi_t(p)\in U$ for sufficiently small $t>0$. 
We denote a local coordinate system  on $U$ by $(x, y)$ and we put the coordinate of $p$ and  $\varphi_t(p)$ by $(x, y)$ and  $(\overline{x}, \overline{y})$, respectively.

Since we have 
\begin{eqnarray*}
  F(\overline{x}, \overline{y})&=&F(x, y)+t \frac{d F(\overline{x}, \overline{y})}{dt}\bigg|_{t=0}+O(t^2)\\
                        &=&F(x, y)+t\bigg(\frac{\partial F(\overline{x}, \overline{y})}{\partial \overline{x}^s}\frac{ d\overline{x}^s}{dt}\bigg)_{t=0}
                                                            +t\bigg(\frac{\partial F(\overline{x}, \overline{y})}{\partial \overline{y}^s}\frac{ d\overline{y}^s}{dt}\bigg)_{t=0}+O(t^2)\nonumber\\
                           &=&F(x, y)+t\frac{\partial F(x, y)}{\partial x^s}W^s +t\frac{\partial F(x, y)}{\partial y^s}\frac{\partial W^s}{\partial x^u}y^u+O(t^2), \nonumber
\end{eqnarray*}
we have

\begin{lemma}\label{Lemma 4.4}
Let $(M, F=\alpha^2/\beta)$ be a  Kropina space whose navigation data is $(h, W)$.  And let $\{\varphi_t\}_{t\in I_\epsilon}$ be a local  $1$-parameter group  of local transformations on $M$ which generate the vector field $W$. 
Then  $\varphi_t$ $(\epsilon>t>0)$  are Finslerian local  isometries of a Kropina metric $F=\alpha^2/\beta$ if and only if $W$ satisfies the differential equation (\cite{LCM})
\begin{eqnarray}\label{equation 4.12}
K_W(F):=\frac{\partial F(x, y)}{\partial x^s}W^s +\frac{\partial F(x, y)}{\partial y^s}\frac{\partial W^s}{\partial x^u}y^u=0.
\end{eqnarray}
\end{lemma}

\begin{definition}(\cite{LCM})
We call $K_W(F)=0$ a Killing equation of a Finsler metric $F$.
\end{definition}

In general, let $F$ be an $(\alpha, \beta)$-metric with defining by  $F=\alpha \phi(s)$, $s=\beta/\alpha$, where $\phi(s)$ is a smooth function of $s$.
Suppose that $\varphi_t$ $(\epsilon>t>0)$  are Finslerian local isometries on a Finsler space $(M, F)$ which satisfies 
\begin{eqnarray*}
  \overline{x}^i =x^i+tV^i+O(t^2),  \hspace{0.1in}
 \overline{y}^i =y^i+t\frac{\partial V^i}{\partial x^j}y^j+O(t^2),
\end{eqnarray*}
then in the similar way to the case of Kropina metric it follows that
  $\varphi_t$ $(\epsilon>t>0)$  are Finslerian local isometries on $(M, F)$ if and only if $V$ satisfies the differential equation
\begin{eqnarray*}
K_V(F):=\frac{\partial F(x, y)}{\partial x^s}V^s +\frac{\partial F(x, y)}{\partial y^s}\frac{\partial V^s}{\partial x^u}y^u=0.
\end{eqnarray*}

Then, the Killing equation $K_V(\alpha \phi(s) )$ reads 
\begin{eqnarray}\label{equation 4.13}
   \bigg(\phi(s)-s\frac{d \phi}{ds} \bigg)K_V(\alpha)+\frac{d \phi}{ds}K_V(\beta)=0.
\end{eqnarray}

Hence, we get
\begin{lemma}(\cite{LCM})\label{lemma 4.5}
For an $(\alpha, \beta)$-metric space $(M, F=\alpha \phi(s))$, where $s=\beta/\alpha$, the Killing equation $K_V(F)$ can be represented by the equation (\ref{equation 4.13}).
\end{lemma}

And  we obtain
\begin{lemma}(\cite{LCM})\label{lemma 4.6}
\begin{eqnarray*}
K_V(\alpha)=\frac{1}{2\alpha}(V_{i;j}+V_{j;i})y^iy^j,\hspace{0.2in}
K_V(\beta)=\bigg(b_{j;i}V^i+b^i V_{i;j}\bigg)y^j,
\end{eqnarray*}
where $";"$ denotes the covariant derivative with respect to the Riemannian metric $\alpha$ and $V_i:=a_{ij}V^j$.
\end{lemma}

Therefore, Lemma \ref{Lemma 4.4} can be written as follows:
\begin{lemma}\label{Lemma 4.7}
Let $(M, F=\alpha^2/\beta)$ be a  Kropina space whose navigation data is $(h, W)$. And let $\{\varphi_t\}_{t\in I_\epsilon}$  be a local  $1$-parameter group  of local transformations on $M$ which generate the vector field $W$. 

Then  $\varphi_t$ $(\epsilon>t>0)$  are Finslerian local isometries of a Kropina metric $F=\alpha^2/\beta$ if and only if $W$ satisfies the differential equation
\begin{eqnarray}\label{equation 4.14} 
2\beta K_W(\alpha)-\alpha K_W(\beta)=0.
\end{eqnarray}
\end{lemma}

 From (\ref{equation 1.1}), we have
\begin{eqnarray*}
 V^i=W^i=\frac{1}{2}b^i, \hspace{0.1in}  V_i=\frac{1}{2}a_{is}b^s=\frac{1}{2}b_i.
\end{eqnarray*}
Using the equations in Lemma \ref{lemma 4.6}, we can change the equation (\ref{equation 4.14}) as follows:
\begin{eqnarray}
\beta \bigg(b_{t;s}+b_{t;s}\bigg)y^sy^t-\alpha^2 \bigg(b_{t;s}b^s+b^s b_{s;t}\bigg)y^t=0.
\end{eqnarray}
Since $\alpha^2$ is not divisible by $\beta$,  the two equations
\begin{eqnarray*}
\bigg(b_{t;s}b^s+b^s b_{s;t}\bigg)y^t=c(x)\beta,    \hspace{0.1in} and  \hspace{0.1in}
\bigg(b_{t;s}+b_{s;t}\bigg)y^sy^t=c(x)\alpha^2           
\end{eqnarray*}
hold for a function $c(x)$ of $(x^i)$ alone.
From the above  equations, we get
\begin{eqnarray} \label{equation 4.16}
b_{t;s}b^s+b^s b_{s;t}=c(x)b_t,          \hspace{0.1in} and  \hspace{0.1in}
b_{t;s}+b_{s;t}=c(x)a_{st}. 
\end{eqnarray}
The first equation is gotten from the second equation. Therefore, the vector field $b^\#= b^i(\partial/\partial x^i)$ satisfies the second  equation  in (\ref{equation 4.16}).

Conversely, suppose that the vector field $b^\#= b^i(\partial/\partial x^i)$ satisfies the second equation  in (\ref{equation 4.16}) holds, then we get the equation (\ref{equation 4.14}).

Hence, we get

\begin{proposition}\label{Theorem 4.8}
Let $(M, F=\alpha^2/\beta)$ be a  Kropina space whose navigation data is $(h, W)$.  Suppose that $\{\varphi_t\}_{t\in I_\epsilon}$ is a local  $1$-parameter group  of local  transformations on $M$ which generates the vector field $W$. 
Then  $\varphi_t$ $(\epsilon>t>0)$ are Finslerian  local isometries of a Kropina metric $F=\alpha^2/\beta$ if and only if the vector field $b^\#= b^i(\partial/\partial x^i)$ satisfies  the second equation in (\ref{equation 4.16}).
\end{proposition}

From Theorem \ref{Theorem 2.1}, Theorem \ref{Theorem 2.2} and Proposition \ref{Theorem 4.8}, we obtain

\begin{theorem}\label{Corollary 4.9}
Let $(M, F=\alpha^2/\beta)$ be a  Kropina space whose navigation data is $(h, W)$.  Suppose that $\{\varphi_t\}_{t\in I_\epsilon}$ is a  local $1$-parameter group   of local transformations on $M$ which generate the vector field $W$. 

 Then, $\varphi_t$ $(\epsilon>t>0)$ are Finslerian local isometries of a Kropina metric $F=\alpha^2/\beta$ if and only if the Kropina space $(M, F)$ is a strong Kropina one.

Namely, $\varphi_t$ $(\epsilon>t>0)$ is a local $1$-parameter group of  Finslerian local isometries of a Kropina metric $F=\alpha^2/\beta$ if and only if  $\varphi_t$ is a local $1$-parameter group of  local isometries on $(M, h)$.
\end{theorem}

\begin{remark}
If $M$ is a 3-dimensional compact connected, simply connected manifold, or a compact Lie group of any dimension, then all isometries of the strong Kropina space $(M,F)$ are global. 

Indeed,
let us remark that the local isometries of the Kropina metric are global when the Killing vector field is complete. This happens for example when $M$ is compact (Proposition \ref{Prop. compact}). 

On the other hand we have shown in \cite{YS} that there are restrictions on the existence of strong Kropina structures on compact manifolds. 

\end{remark}


\vspace{0.5in}

\begin{center}
          Hachiman technical\\         
      High School \\                 
	5 Nishinosho-cho Hachiman\\     
      523-0816 Japan \\              
	E-mail: ryozo@e-omi.ne.jp
  \end{center}	
\begin{center}
     Department of Mathematics\\
     Tokai University\\
       Sapporo, 005\,--\,8601 Japan\\
        E-mail:    sorin@tspirit.tokai-u.jp
\end{center}	
\end{document}